\documentclass[12pt]{amsart}
\usepackage[margin=1in]{geometry}                % See geometry.pdf to learn the layout options. There are lots.
\usepackage{graphicx}
\usepackage{comment}

\usepackage{amsmath}
\usepackage{amssymb}
\usepackage{amsfonts}
\usepackage{amsthm}
\usepackage{mathrsfs}
\usepackage{enumerate}
\usepackage{float}
\usepackage{fancyhdr}
\usepackage{xcolor}
\usepackage{hyperref}
\usepackage{mathtools}
\usepackage{wasysym}

\newtheorem{theorem}{Theorem}
\newtheorem{proposition}[theorem]{Proposition}

\newtheorem{question}[theorem]{Question}

\numberwithin{equation}{section}

\newcommand{\RR}{{\mathbb R}}

\newcommand{\R}{\mathbb{R}}
\newcommand{\Z}{\mathbb{Z}}

\begin{document}
\title{Spectrality of the product does not imply that the components are spectral}
\author{G\'abor Somlai}

\address{Department of Mathematics \\ Eötvös Lor\'and University Faculty of Science, Pázmány Péter sétány 1/C, \\ Budapest, Hungary, H-1117 \\ E-mail: {\tt gabor.somlai@ttk.elte.hu}}
\date{\today}
\begin{abstract}
Greenfeld and Lev conjectured that the Cartesian product of two sets $A$ and $B$ is spectral if and only if $A$ and $B$ are spectral. We construct a counterexample to this conjecture using the existence of a tile that has no spectra. 
\end{abstract}
\maketitle

\section{Introduction}
Let $G$ denote a locally compact abelian group, endowed with Haar measure. A bounded measurable set $S$ (of non-zero measure) is called spectral in $G$, if there exists a collection of characters of $G$ whose  restrictions to $S$ form an orthogonal basis of the function space $L^2(S)$. A bounded measurable set $T$ (of non-zero measure) is called a tile in $G$, if there exist a subset $U\subset G$ such that almost every element of $G$ can  be written uniquely as $t+u$, where $t \in T$ and $u \in U$. 

Both definitions are simpler if $G$ is a finite abelian group. In this case the dual group $\hat{G}$ is isomorphic to $G$. If a subset of a finite abelian group is spectral, then its spectrum is a subset of $\hat{G}$, and it can be identified with a subset of $G$. This allows us to talk about a spectral pair $(A,B)$, where $A \subset G$ and $B\subset \hat{G} \cong G$. Clearly, in this case $|A|=|B|$.
The measure is the counting measure so the only set of measure zero is the empty set.    

It is not too difficult to see that if $A$ and $B$ are both spectral in $G$ and $H$, respectively, then $A \times B$ is spectral in $G \times H$.  
Greenfeld and Lev have shown \cite{GL1} that if $A$ is an interval $\mathbb{R}$ and $B$ is a measurable subset of $\mathbb{R}^{n-1}$ and $A \times B$ is spectral in $\mathbb{R}^n$, then $A$ and $B$ are both spectral. They could use it to prove that for cylindrical convex bodies in $\mathbb{R}^3$ Fuglede's conjecture \cite{F} holds, i.e. spectral sets and tiles coincide.

 Kolountzakis \cite{K1} extended the result of Greenfeld and Lev by proving that if $A$ is the union of two intervals, then the spectrality of $A \times B$ implies the spectrality of $A$ and $B$. Later, Greenfeld and Lev \cite{GL2} proved that the product of a convex polygon $A \subset \R^2$ with a bounded measurable set $B$ is spectral if and only if both $A$ and $B$ are spectral. Finally Kolountzakis, Lev and Matolcsi \cite{KLM} extended this result to the case when $A$ is a convex body. More precisely they proved that if $A$ is a convex body in $\R^n$ and $B$ is any bounded, measurable set in $R^m$ and  $A\times B$ is a spectral set then $A$ must be spectral. Finally, it was proved in \cite{CLZ} that the product of a classical Sierpinski self-affine tile with a Lebesgue measurable set of measure 1 is spectral if and onyl if both components are spectral. 

It is also important to mention that the equivalence of $A \times B$ is a tile and $A$ and $B$ are tiles was observed by Kolountzakis \cite{K1}. 
Kolountzakis asked the next question \cite{K1}, and later repeated \cite{K2} the same question, saying that a negative answer would break the symmetry between the spectral and the tile concepts\footnote{which we think is broken enough by counterexamples to both directions of Fuglede's conjecture. This view is not shared by some of the more experienced people in the field.}. 
\begin{question}\label{question}
    Let $ \Omega = A \times B$ be the Cartesian product of two bounded, measurable subsets of $A$ of $\RR^n$ and  $B$ of $\RR^m$. When is $\Omega$ spectral? 
\end{question}

Greenfeld and Lev conjectured \cite{GL2} that $\Omega$ is spectral if and only if both $A$ and $B$ are spectral. This conjecture was reiterated in \cite{KLM}\footnote{Question \ref{question} was mentioned in the open problem session of Harmonic and Spectral Analysis conference (2021) by Nir Lev: 
\href{http://mathspectral.hu/wp-content/uploads/2021/06/problems_and_remarks_HSA2021.pdf}{HSA2021 problems}}.

We disprove this conjecture by first constructing a counterexample for finite abelian groups and then lifting it up to Euclidean space along the lines of Tao's idea \cite{T}, using the standard technology developed by Kolountzakis and Matolcsi \cite{KM}. 
\begin{theorem}\label{thm:disprove}
There is a pair of bounded measurable sets $A$ and $B$ in $\RR^3$ such that $A$ is not spectral but $A \times B$ is spectral in $\RR^6$.  
\end{theorem}
In order to disprove this conjecture we prove the following theorem for finite abelian groups. This the main ingredient of this paper since Theorem \ref{thm:disprove} follows using standard methods we will recall. 

\begin{theorem}\label{thm:finite tile}
Let $G$ be a finite abelian group of order $n$ and let
$D=\left\{ (g,g) \mid g\in G \right\}$ denote the diagonal subgroup of $G \times G$.
\begin{enumerate}
    \item\label{item1}  Let $P=\left\{ (a_i,b_i) \mid i=1, \ldots, n \right\}$ be a subset of $G \times G$. Then $(P,D)$ is a spectral pair if and only if $\left\{ a_i+b_i \mid i=1, \ldots, n \right\}=G$.
\item\label{item2} Let $A$ and $B$ subsets of $G$ with $|A|*|B|=|G|$. Then $A$ tiles with $B$ if and only if $(A \times B,D)$ is spectral in $G \times G$.
\end{enumerate}
\end{theorem}
This theorem tells us that besides cardinality considerations tiling partners in a group $G$ can be described by the spectrality of one single set with a fixed spectrum, which is the diagonal of $G \times G$.
\section{Proof of the results}
 In order to verify that a subset $S$ of $G\times G$ is spectral, we need to find a set of characters of $G \times G$ such that their restrictions to $S$ form a basis in the space of complex valued functions on $S$. Since the dimension of this space is $|S|$ it is sufficient to find $|S|$ exponential functions whose restriction to $S$ is pairwise orthogonal. 

Now we prove Theorem \ref{thm:finite tile}. 
\begin{enumerate}
    \item Let us first assume that $\left\{ a_i+b_i \mid i=1, \ldots, n \right\}=G$. It is clear that both $P$ and 
$D=\left\{ (g,g) \mid g\in G \right\}$  are of cardinality $|G|$ so it is enough to verify that the representations corresponding  to the elements of $D$ are pairwise orthogonal on $A \times B$. Let $0 \ne g \in G$. The following computation shows the required orthogonality of the characters if $\left\{ a_i+b_i \mid i=1, \ldots, n \right\}=G$. 
\begin{equation}\label{eq:21}
    \begin{split}
 &  \chi_{(g,g)}(P) = \sum_{i=1}^n  \chi_{(g,g)}(a_i,b_i)= 
   \sum_{i=1}^n \chi_g(a_i) \chi_g(b_i)= \\
  & \sum_{i=1}^n \chi_g(a_i+b_i)  = \sum_{h \in G} \chi_g(h)=0 . 
    \end{split}
\end{equation}
%The converse of the statement also follows easily from this computation since a function on $G$ is constant if and only if the Fourier transformation of the function vanishes on every nontrivial character. 
Let us assume that $(P,D)$ is a spectral pair. Then by equation \eqref{eq:21} we have that $\sum_{i=1}^n \chi_g(a_i+b_i)  =0$ for every $0 \ne g \in G$. Thus the multiset $\left\{ a_i+b_i \mid i=1,\ldots, n \right\}$ is constant on $G$. Since $G$ has exactly $n$ elements we have that  $\left\{ a_i+b_i \mid i=1,\ldots, n \right\}=G$. 
%\textcolor{red}{Ezt jobban el kene magyarazni, hogy a forditott iranyu kovetkeztetes miert igaz. Ez a mondat jelenleg elegge "confusing", hogy milyen konstans fuggvenyrol beszelsz.}
\item 
This statement directly follows from Theorem \ref{thm:finite tile} \eqref{item1}. \qed
\end{enumerate}
Now we prove Theorem \ref{thm:disprove}.

The first example of a non-spectral tile was found by Kolountzakis and Matolcsi \cite{KM}. The complexity of the finite abelian groups, where tiles with no spectra can be found was reduced in \cite{FMM} and \cite{FR}. Counterexamples for the tile-spectral direction of Fuglede's conjecture in $(\mathbb{Z}_{24})^3$, which is the direct sum of 3 cyclic groups, was found in \cite{FMM}. 

The existence of such an example (tile with no spectra) combined with Theorem \ref{thm:finite tile} \eqref{item2} provides a pair of subsets of $(\mathbb{Z}_{24})^3$, one of which is non-spectral such that the Cartesian product of the two sets is spectral. Let $A$ be a non-spectral tile in $(\mathbb{Z}_{24})^3$ with one of its tiling complement $B$.

We now lift the construction from finite abelian groups to a counterexample to the original conjecture in $\mathbb{Z}^6$ and $\mathbb{R}^6$.
We borrow the notation of \cite{KM} so let us assume that $A$ is a subset of $\{0,1, \ldots, n_1-1\}\times \ldots \times \{0,1, \ldots, n_d-1\}$. Then for a positive integer $k$ we denote by $A(k)$ the sum of $A$
and $$\{0,n_1,2n_1,\ldots,(k-1)n_1\} \times\ldots \times \{0,n_d,2n_d,\ldots,(k-1)n_d\}.$$

Let $A$ be a subset of $\Z_{24}^3$, which is a tile and not spectral. We may naturally identify $A$ with a subset of $\{0,1,\ldots, 23 \}^3 \subset \mathbb{Z}^3$.
Let $B$ one of $A$'s tiling partners in $\Z_{24}^3$, which we identify with a subset of $\{0,1,\ldots, 23 \}^3 \subset \mathbb{Z}^3$. We may also embed $A\times B$ into $\mathbb{Z}^6$.

It is easy to see that for every positive integer $k$ we have $A(k) \times B(k)= (A\times B)(k)$. 
One can also verify that if $C$ is a spectral set in $\mathbb{Z}_{24}^6$, then if we embed $C$ into $\Z^6$ as above we obtain a spectral set with a rational spectrum.
Since $A\times B$ is spectral in $\Z_{24}^6$, we have that $A \times B$ is spectral in $\Z^6$. It follows from Proposition 2.1 in \cite{M} that $(A\times B)(k)$ is also spectral in $\mathbb{Z}^6$, for every $k$. 

On the other hand, Theorem 4.1 in \cite{KM} shows that 
$A(k)$ is not spectral in $\mathbb{Z}^6$ if $k$ is large enough, since $A$ is not spectral in $(\mathbb{Z}_{24})^3$. Now we may conclude that there is a pair of finite subsets of $\Z^3$, one of them is not spectral, such that the product of the two sets is spectral in $\Z^6$.

Theorem 4.2 in \cite{KM} shows that $A(k)+[0,1)^3$ is not spectral if $k$ is large enough while $(A \times B)(k) +[0,1)^6$ is spectral in $\R^6$ for every $k$. Furthermore,  ($A(k)+[0,1)^3) \times (B(k)+[0,1)^3)=(A \times B)(k) +[0,1)^6$, finishing the proof of Theorem \ref{thm:disprove}. \qed

Finally, one can verify that $G=\left\{a_i+b_i \mid 1 \le i \le n \right\}$
holds if and only of $P$ is a coset representative of $\left\{ (g,-g) \mid g \in G \right\}$ so it tiles with a subgroup.
This shows that the spectra of the diagonal $D$ in $G \times G$ have a common tiling complement. Thus this new example of a spectrum does not (directly) lead to a counterexample for the spectral-tile direction of Fuglede's conjecture.

%\noindent{{\sc G\'abor Somlai:}  \\
%ELTE-TTK, Institute of Mathematics \\ P\'azm\'any P\'eter s\'et\'any 1/C, Budapest, Hungary, H-1117},\\
%E-mail: {\tt gabor.somlai@ttk.elte.hu}
\section*{Acknowledgements}
Research supported  by the Hungarian National Foundation for
Scientific Research, Grant: 138596.

The author is grateful to M\'at\'e Matolcsi for his help in transforming the first draft/version of the manuscript into a meaningful manuscript.

The author is grateful to Nir Lev for his suggestions on the presentation and mathematical content of the manuscript.  

The author is grateful to Azita Mayeli, who brought this problem to his attention again a year ago. 
\end{document}